\documentclass[fleqn]{article}

\title{On Poincar\'{e}'s Fourth and Fifth Examples of Limit Cycles at Infinity}
\author{Roland K. W. Roeder \\
	University of California, San Diego}

\date{January 14, 2001}

\input{psfig.sty}

\begin{document} \maketitle \begin{abstract} Errors are found in example
problems from Henri Poincar\'e's paper ``M\'{e}moire sur les courbes
d\'{e}finies par une \'{e}quation diff\'{e}rentielle.'' Examples four and
five from chapter seven and examples one, two, and three from chapter nine
do not have the limit cycles at infinity predicted by Poincar\'e.  
Instead they have fixed points at every point at infinity.  In order to
understand the errors made by Poincar\'e, examples four and five are
studied at length. Replacement equations for the fourth and fifth examples
are suggested based on the supposition that terms were omitted from
Poincar\'e's equations. \end{abstract}

\section{Introduction} In the late nineteenth and early twentieth century
Henri Poincar\'{e} began to study the qualitative aspects of systems of
differential equations.  This analysis was a breakthrough in the field
because one no longer had to obtain a specific solution of the equation in
question to understand its general behavior.  This manner of analysis was
introduced by Poincar\'{e} in his paper ``M\'{e}moire sur les courbes
d\'{e}finies par une \'{e}quation diff\'{e}rentielle'' \cite[1881]{PON}.  
Poincar\'{e} worked primarily with systems of two variables and played a
key role in identifying the existence of limit cycles with the
Poincar\'{e}-Bendixson theorem.  Poincar\'{e} also searched for a complete
global analysis of a system of two variables; to do so he introduced
analysis at infinity by means of the Poincar\'{e} Sphere.  These two
aspects of his analysis join in a behavior known as a limit cycle at
infinity.

To illustrate the application of these techniques, Poincar\'e presents
example problems in chapters seven and nine. In chapter seven,
Poincar\'{e} presents five example problems. He states that examples
three, four, and five have limit cycles at infinity.  (Examples three,
four and five throughout the paper will be referred to as P3, P4, and P5
respectively.) However, when one checks Poincar\'{e}'s assertion, one
finds that examples P4 and P5, in fact, do not have limit cycles at
infinity.  The errors remains in the collected works \cite{PON2} and the
error in P4 remains in a modern work that includes translations of
examples P3 and P4. \cite[p. 155]{DAV}.

In chapter nine, Poincar\'e presents three additional example problems.  
He states that each of these examples have limit cycles at infinity.  
However, he makes the same error analyzing these equations as he did for
P4 and P5.

Throughout this paper, in an effort to understand why Poincar\'e made these
errors, we focus on studying examples P4 and P5 from chapter seven.
In the conclusion we will attempt to apply what we have learned about 
the errors in P4 and P5 to the examples from chapter nine.

Because the analysis at infinity is algebraically simple for each of these
problems, one becomes perplexed at this finding. Surely Poincar\'{e} could
not have made such an error. With this belief, one might suppose that
there was a simple omission of some sort.  If one were to find equations
that have the same behavior as Poincar\'{e} predicts for P4 and P5, and
that, were one to drop certain terms, could become his printed equations,
this would serve as a plausible explanation.  We look at P4 in the context
of P3, and consider a plausible modification. Such an analysis is
encouraged by the strong geometric similarities between the two systems.  
The modification restores the limit cycle at infinity, and hence supports
the hypothesis of an error of omission, becoming a key candidate for the
intended fourth example. We will call this example R4.  Example P5 does
not fit in the same genre of geometric system and hence no such
modification suggests itself.

When one investigates the uniqueness of these ``solutions'' to the problem
of finding modifications of Poincar\'{e}'s equation that match his
analysis for P4, one can find a classification of all modifications of the
original equation, of a certain type, that result in the correct behavior.  
There exist a great number of such modifications;  however, this
remarkable non-uniqueness in ``solutions'' does not cause skepticism that
R4 is the intended equation for the fourth example.  It is seen that the
modification to obtain R4 is the simplest of the modifications found under
the classification, and hence, through application of the philosophical
principle known as Occam's razor, R4 remains the most plausible equation
for the fourth example.

Based on the success of the technique applied to P4, a similar technique
is applied to P5 to try to determine a possible modification of P5 which
yields the correct behavior.  This technique provides a number of
modifications which agree with Poincar\'{e}'s analysis for P5 thus
illustrating the power of the technique.  One of these, which we shall
call R5, is chosen as a suggested replacement to P5.  (No modifications
are suggested for the examples from chapter nine.)

To aid the reader with the analysis of differential equations at infinity,
a section is included presenting the standard techniques for analysis
of the  behavior at infinity for a polynomial system in rectangular 
coordinates.  A more geometric analysis is also presented that is 
quite valuable in many cases in which such a system simplifies in polar 
coordinates.
     
\setcounter{equation}{0} \section{Differential Equations at Infinity}
Differential equations that are defined on the plane can be analyzed at
infinity by extending the plane to a representation of the plane that
includes the ``points at infinity.'' This analysis began with the work of
Poincar\'{e} and Bendixson \cite{PON}\cite{BEN}. Throughout this paper we
use Poincar\'e's method of analysis.  More modern and readable accounts on
the Poincar\'e's method are available in books by Perko \cite{PER},
Lefschetz \cite{LEF}, and Minorsky \cite{MIN}.  This section begins with a
brief summary of the standard results as in Perko \cite[pp. 264-268]{PER}.
A theorem is presented, consistent with the standard techniques, based on
polar coordinates, that demonstrates the geometric nature of these systems
at infinity.  This theorem often makes it possible for one to determine by
inspection whether a system has a limit cycle at infinity.

Through the use of projective geometry, the complete behavior of a two
dimensional system of differential equations can be seen in its behavior
on a sphere of finite radius known as the Poincar\'{e} sphere. To do this,
one places the phase plane tangent to the sphere and makes correspond
points on the plane with the points on the sphere by central projection (a
point on the plane corresponds with an antipodal pair on the sphere.)  
One must note that the line intersects the sphere at two points; to remove
this non-uniqueness, antipodal points are identified on the Poincar\'{e}
sphere.  The points that were at infinity on the original plane become
points on the equator of the sphere.

From the behavior on the Poincar\'{e} sphere, one can construct the
``global phase portrait.'' To consider the global phase portrait of a
system, one projects the trajectories from the upper hemisphere
orthogonally down onto the plane that goes through the center of the
sphere and is parallel to the original plane.  In this way the complete
behavior on the original plane becomes the behavior on the finite disk of
this new plane.  Points at infinity become the points at which the sphere
intersects the plane -- the boundary of the disk.

The first analysis that one would do for a system in the finite plane is
to find the location of the fixed points.  With this as motivation, the
primary problem of analyzing the behavior of a system at infinity is to
determine the location of fixed points, if there are any.  One considers
the system:  
\begin{eqnarray*} 
\frac{dx}{dt} & = & P(x,y), \\
\frac{dy}{dt} & = & Q(x,y), 
\end{eqnarray*} 
where $P$ and $Q$ are polynomials of degree $m$ in $x$ and $y$.  Denote by
$P_m$ and $Q_m$ the homogeneous polynomials consisting of terms of degree
exactly $m$.
 
\newtheorem{theorem}{Theorem} 
\begin{theorem} 
The critical points at infinity for the system above occur at the points
(X,Y,0) on the equator of the Poincar\'{e} Sphere where $X^2+Y^2=1$ and
\begin{eqnarray} 
XQ_m(X,Y)-YP_m(X,Y)=0 \label{PER} 
\end{eqnarray} 
or equivalently at the polar angles $\theta$ and $\theta + \pi$ which are
solutions of 
\begin{eqnarray} 
G_{m+1}(\theta) \equiv cos\theta \,
Q_m(cos\theta ,sin\theta )-sin\theta \, P_m(cos\theta ,sin\theta )=0.  
\label{G} 
\end{eqnarray}
This equation has at most $m+1$ pairs of roots $\theta$ and $\theta +
\pi$ unless $G_{m+1}(\theta )$ is identically zero.  If $G_{m+1}(\theta )$
is not identically zero, then the flow on the equator of the Poincar\'{e}
sphere is counter-clockwise at points corresponding to polar angles
$\theta $ where $G_{m+1}(\theta ) > 0$ and it is clockwise at points
corresponding to polar angles where $G_{m+1}(\theta ) < 0$. 
\end{theorem}
See Perko \cite{PER}, section 3.10, Theorem 1.

An important consequence of this Theorem 1 is that any polynomial system
in rectangular coordinates can be extended to the Poincar\'e sphere.

One should be aware that in the proof of Theorem 1, a new time scale is
defined to study the trajectories in the proximity of the equator.  There
are two opposing conventions used in the degenerate case where (\ref{PER})
is identically zero. The most common one, used by Perko \cite{PER} and
Minorsky \cite{MIN} is to re-parameterize time in such a way that the
equator is forced to consist of trajectories and fixed points.  The second
convention, followed by Lefschetz \cite{LEF} and Poincar\'e \cite{PON}
uses a different time scale, allowing trajectories to cross the equator
(implying that trajectories in one of the hemispheres will flow in a
direction inconsistent with the flow on the plane.) Because the former
method is more common, we will adhere to it, even though Poincar\'e used
the latter.  (A brief comment will be made in the analysis of Example P4
about why this does not invalidate our results.)

If there are no fixed points at infinity there is a cycle at infinity.  
If trajectories in the proximity of this cycle at infinity approach (or
recede from) it, it is a limit cycle at infinity. A corollary immediately
follows:

\newtheorem{corollary}{Corollary} 
\begin{corollary} 
Let $r$ denote radial distance in polar coordinates. If the system
considered above has no fixed points at infinity and $dr/dt \neq 0$ for $r
\geq R$, for some $R$, then there exists a limit cycle at infinity.  
\end{corollary}

While this theorem makes it relatively easy to determine the behavior of a
system at infinity, there is a more geometric analysis that for certain
systems will make the behavior at infinity particularly easy to determine.

\begin{theorem}
Consider a polynomial system
\begin{eqnarray}
\hspace{1.0in} \frac{dx}{dt}=P(x,y), \hspace{.25in} \frac{dy}{dt}=Q(x,y) \label{PQ}
\end{eqnarray}
with expressions for $dr/dt$ and $d\theta/dt$ of order $I$ and $J$ in $r$
as $r \rightarrow \infty$, respectively.  Let $k=I-J$. Then, if $k\geq2$,
the equator of the Poincar\'{e} sphere consists entirely of fixed points.
If $k\leq 1$ then: if $G_{m+1}(\theta)\neq 0$ for all $\theta$ there is a
cycle at infinity.  Furthermore, if $dr/dt$ satisfies the conditions of
the above corollary, it is a limit cycle. Otherwise, the equator of the
Poincar\'{e} sphere has finitely many fixed points located at $\theta$
such that $G_{m+1}(\theta)=0$, as in Theorem 1.
\end{theorem}
Proof of Theorem 2: \newline Suppose that $P$ and $Q$ are polynomials of
degree $m$ in $x$ and $y$.  We can express this system in polar
coordinates, for $r\neq 0$, as follows, 
$$\frac{dr}{dt}=(cos\theta P_0+sin\theta Q_0)+r(cos\theta P_1(cos\theta,
  sin\theta)+sin\theta Q_1(cos\theta,sin\theta))+ $$
$$\cdots+r^m(cos\theta P_m(cos\theta,sin\theta)+sin\theta Q_m(cos\theta,
  sin\theta)). $$
$$\frac{d\theta}{dt}=r^{-1}(cos\theta Q_0-sin\theta P_0)+(cos\theta Q_1(cos
  \theta,sin\theta)-sin\theta P_1(cos\theta,sin\theta))+ $$
$$\cdots+r^{m-1}(cos\theta Q_m(cos\theta,sin\theta)-sin\theta P_m(cos
  \theta,sin\theta)).$$
Throughout the remainder of the proof we require $r\neq0$.
The following definitions will greatly simplify notation:
\begin{eqnarray*}
\eta_0(\theta) & = & (cos\theta P_0+sin\theta Q_0), \\
\hspace{1in} \vdots  \\
\eta_I(\theta) & = & (cos\theta P_I(\cos\theta,sin\theta)+sin\theta
		Q_I(cos\theta,sin\theta)),
\end{eqnarray*}
and
\begin{eqnarray*}
\xi_{-1}(\theta) & = & (cos\theta Q_0-sin\theta P_0), \\
\hspace{1in} \vdots \\
\xi_J(\theta) & = & (cos\theta Q_{J+1}(\cos\theta,sin\theta)-sin\theta
		P_{J+1}(cos\theta,sin\theta)),
\end{eqnarray*}
where $I$ and $J$ are the highest degree terms, in $r$, in the equations
for $dr/dt$ and $d\theta/dt$ respectively. With this notation we obtain,
\begin{eqnarray*}
\frac{dr}{dt} & = & \eta_0(\theta)+\cdots+r^I \eta_I(\theta), \\
\frac{d\theta}{dt} & = & r^{-1}\xi_{-1}(\theta)+\xi_0(\theta)+\cdots+r^J \xi_J(\theta).
\end{eqnarray*}
By the definition of $k$ we have $k=I-J$. If we put this system in
differential form we obtain 
\begin{eqnarray}
(r^{-1}\xi_{-1}(\theta)+\cdots+r^J \xi_J(\theta))dr - 
(\eta_0(\theta)+\cdots+ r^I \eta_I(\theta))d\theta = 0. \label{DF1}
\end{eqnarray}

To determine the behavior at infinity, one projects a differential
equation onto the Poincar\'{e} sphere.  To understand the details of the
behavior at infinity, project the upper hemisphere of the Poincar\'e
sphere onto the cylinder of radius $1$ with axis orthogonal to the phase
plane that is tangent to the Poincar\'{e} sphere at the equator.  (We can
restrict our attention to the upper hemisphere because we follow the
standard convention, which does not allow trajectories to cross the
equator.)  The ``equator'' of the cylinder is the part of the cylinder
that touches the sphere.  By geometric analysis of the projection from the
plane directly to the cylinder, $s=1/r$ and correspondingly
$dr=\frac{-1}{s^2}ds$.  This projection leaves $\theta$ unchanged.

Projecting (\ref{DF1}) from the plane directly onto the cylinder gives:
$$(s \xi_{-1}(\theta)+\cdots+\frac{1}{s^J} \xi_J(\theta))(\frac{-1}{s^2}ds) - 
  (\eta_0(\theta)+\cdots+\frac{1}{s^I} \eta_I(\theta))d\theta = 0,$$
which simplifies to
\begin{eqnarray}
(s \xi_{-1}(\theta)+\cdots+\frac{1}{s^J} \xi_J(\theta))ds + 
 s^2(\eta_0(\theta)+\cdots+\frac{1}{s^I} \eta_I(\theta))d\theta = 0. \label{DF2}
\end{eqnarray}

\noindent{\bf{Case 1:} $k\geq 2$}

In this case we have $I-J\geq 2$, so $I-2 \geq J$.  We multiply (\ref{DF2})
by $s^{I-2}$ to clear the denominators obtaining:
\begin{eqnarray}
(s^{I-1} \xi_{-1}(\theta)+\cdots+s^{I-J-2} \xi_J(\theta))ds +
  (s^I \eta_0(\theta)+\cdots+ \eta_I(\theta))d\theta = 0 \label{DFC1A}
\end{eqnarray}
However, this equation would indicate trajectories crossing the equator
(from one side or the other) for all points $(\theta,0)$ where
$\eta_I(\theta) \neq 0$.  Because we wish to maintain the same convention
as in Theorem 1, the equator is required to consist only of fixed points
and trajectories moving along the equator. As a result, the trajectories
crossing the equator would lead to a violation of the uniqueness of
solutions to (\ref{DFC1A}) at these points $(\theta,0)$.  To resolve this,
we must multiply (\ref{DFC1A}) by an additional value of $s$, a
reparameterization of time that causes these trajectories to slow down so
that they do not cross the equator.  We obtain:  
\begin{eqnarray} 
(s^I \xi_{-1}(\theta)+\cdots+s^{I-1-J} \xi_J(\theta))ds +
  (s^{I+1} \eta_0(\theta)+\cdots+s \eta_I(\theta))d\theta = 0 \label{DFC1B}
\end{eqnarray}
With regards to this new time scale, $\tau$, (\ref{DFC1B}) can be
expressed as the system:
\begin{eqnarray*}
\frac{ds}{d\tau} & = &-(s^{I+1}\eta_0(\theta)+\cdots+s\eta_I(\theta)), 
\end{eqnarray*}
\begin{eqnarray}
\frac{d\theta}{d\tau}& = & (s^I \xi_{-1}(\theta)+\cdots+s^{I-J-1}\xi_J(\theta)).
\label{DC1}
\end{eqnarray}
Because $k=I-J \geq 2$, at $s=0$, we have:
$$\frac{ds}{d\tau}=0,\hspace{.25in} \frac{d\theta}{d\tau}=0.$$
We conclude that every point on the equator is a fixed point.

\noindent{\bf{Case 2:} $k\leq 1$}

We have $I-J\leq1$ and $J\leq I-1$.  We multiply (\ref{DF2}) by $s^J$ to clear
the denominators.  In this case, no trajectories approach the equator in a 
finite time, hence there is not need for the re-parameterization done in
case 1.

We obtain:
$$(s^{J+1} \xi_{-1}(\theta)+\cdots+s \xi_{J-1}(\theta)+\xi_J(\theta))ds + 
  (s^{J+2} \eta_0(\theta)+\cdots+s^{2-(I-J)}\eta_I(\theta))d\theta = 0.$$
Expressing this in the form of a system, with respect to the new time
scale $\tau$, we obtain: 
\begin{eqnarray*} 
\frac{ds}{d\tau} & = &
-s^{2-(I-J)}(s^I \eta_0(\theta)+\cdots+\eta_I(\theta)), \\
\frac{d\theta}{d\tau} & = & (s^{J+1}
\xi_{-1}(\theta)+\cdots+\xi_J(\theta)). 
\end{eqnarray*} 
So, on the equator we have $s=0$ and $$\frac{ds}{d\tau}=0, \hspace{.25in}
\frac{d\theta}{d\tau}=\xi_J(\theta).$$ So, the fixed points at infinity
are for $\theta$ such that $\xi_J(\theta)=0$. It is easy to see that, in
this case, $\xi_J(\theta)=G_{m+1}(\theta)$, as defined in Theorem 1.  
Hence, the fixed points at infinity are at $\theta$ such that
$$G_{m+1}(\theta)=cos\theta Q_m(cos\theta,sin\theta)-sin\theta
  	P_m(cos\theta,sin\theta)=0.$$
We conclude that if $G_{m+1}(\theta) \neq 0$ for all $\theta$ then there is
a cycle at infinity, and that, if nearby trajectories approach that cycle
(or recede from it), as in the corollary to Theorem 1, then it is a limit 
cycle.  

This concludes a proof of Theorem 2.

One should notice that in the proof of this theorem one finds that for the
case $k \leq 1$ one has $G_{m+1}(\theta)=\xi_J(\theta)$, the highest order
term in $r$ of the $d\theta/dt$ equation.  This makes this theorem
particularly useful -- one can often tell whether a system has a limit
cycle at infinity by merely expressing the system in polar coordinates.

One should also note that in the proof of case 1, if $k > 2$, equation
(\ref{DC1}) gives that a trajectory approaching (or receding from) any
point on the equator does so orthogonally to the equator.

\setcounter{figure}{0}
\section{The Poincar\'{e} Examples}
To enable the reader to understand Poincar\'{e}'s fourth and fifth
examples, three of Poincar\'{e}'s examples from chapter seven of his paper
\cite[pp 274-281]{PON} are presented.  Poincar\'{e}'s third example, P3,
works as an introduction to his fourth and fifth examples and serves as an
example of a system that has a limit cycle at infinity, the behavior that
Poincar\'{e} claims for P4 and P5. Examples P4 and P5 are then presented
with demonstration of the error in his analysis at infinity. At the end of
this section, the error in Poincar\'e's analysis of his examples from
chapter nine is briefly discussed.

\vspace{2 mm}

\noindent{\bf Example P3}

Poincar\'e considers the equation: 

$$\frac{dx}{x(x^2+y^2-1)-y(x^2+y^2+1)} = \frac{dy}{y(x^2+y^2-1)+x(x^2+y^2+1)}.$$
This differential form of the equation is equivalent to the more 
familiar form of a system as noted in \cite[p. 266]{PER}.
\begin{eqnarray*}
\frac{dx}{dt} & = & x(x^2+y^2-1) -y(x^2+y^2+1), \\
\frac{dy}{dt} & = & y(x^2+y^2-1) +x(x^2+y^2+1).
\end{eqnarray*}
The geometric aspects of this system become far more obvious in polar
form:
\begin{eqnarray*}
\frac{dr}{dt} & = & r(r^2-1), \\
\frac{d\theta }{dt} & = & r^2+1.
\end{eqnarray*}
Poincar\'{e} presents the key qualitative features of this system:
\begin{itemize}
\item  The origin: $(x,y)=(0,0)$ is a fixed point.  One has a stable spiral at the origin.
\item   A repelling limit cycle at radius $r=1$ centered about the origin. 
\item   An attracting limit cycle at infinity.
\end{itemize}
All of the features in the finite plane are easily established, so only
the features of this system at infinity are analyzed in detail.  
Poincar\'{e} establishes that there is a limit cycle at infinity because
he claims that there are no fixed points at infinity: ``Il n'y a aucun
point singulier sur l'\'{e}quateur, qui est une caract\'{e}ristique et qui
est par cons\'{e}quent un cycle limite'' \cite[p.~278]{PON}. He omits this
calculation, so we verify it here.

We apply Theorem 2 to show that there are no fixed points at infinity.
From the expression for the P3 in polar coordinates, we find $I=3$ and
$J=2.$ Further, $\xi_2(\theta)=1$ for all $\theta$, demonstrating the
existence of a cycle at infinity. Because $dr/dt > 0$ when $r \geq 2$,
this cycle is a limit cycle.

\vspace{2 mm}

\noindent{\bf Example P4 }

Poincar\'{e}'s fourth example can be considered a more complicated version
of his third example.  Poincar\'{e} writes the equation as:
\begin{eqnarray*} \frac{dx}{x(x^2+y^2-1)(x^2+y^2-9)-y(x^2+y^2-2x-8)} = \\
 \frac{dy}{y(x^2+y^2-1)(x^2+y^2-9)+x(x^2+y^2-2x-8)},
\end{eqnarray*}
which he transforms into the system:
\begin{eqnarray*}
\frac{dx}{dt} & = & x(x^2+y^2-1)(x^2+y^2-9) -y(x^2+y^2-2x-8), \\
\frac{dy}{dt} & = & y(x^2+y^2-1)(x^2+y^2-9) +x(x^2+y^2-2x-8).
\end{eqnarray*}
Poincar\'{e} states the key features of the system: 
\begin{itemize}
\item  The origin: $(x,y)=(0,0)$ is a fixed point.  
One has an unstable spiral at the origin.
\item A: $(x,y)=(1/2,\sqrt{35}/2)$, and B: $(x,y)=(1/2,-\sqrt{35}/2)$ are
fixed points: an unstable node and a saddle respectively, at the
intersection of the circles defined by $x^2+y^2-9=0$ and $x^2+y^2-2x-8=0$.
\item  An attracting limit cycle of radius $r=1$ centered about the origin.
\item  Two heteroclinic orbits connecting fixed points $A$ to fixed point $B$
 along the circle $x^2+y^2-9=0$.
\item  An attracting limit cycle at infinity.
\end{itemize}

All of the stated features in the finite plane can be easily established.
(A similar analysis is applied to a slightly more complicated system in
the following section number five.) The geometric nature of the system
becomes clear in polar coordinates. The system becomes: 
\begin{eqnarray*}
\frac{dr}{dt} & = &r(r^2-1)(r^2-9), \\ 
\frac{d\theta }{dt} & = & r^2-2r\,cos\,\theta - 8. 
\end{eqnarray*} 
The phase portrait corresponding to Poincar\'{e}'s analysis is included in
Figure~3.1. He claims that, as in P3, there are no fixed points at
infinity, and hence a limit cycle. However, when one does the calculations
to check Poincar\'{e}'s assertion, one finds that this is not the case.

\begin{figure}
\centerline{\psfig{file=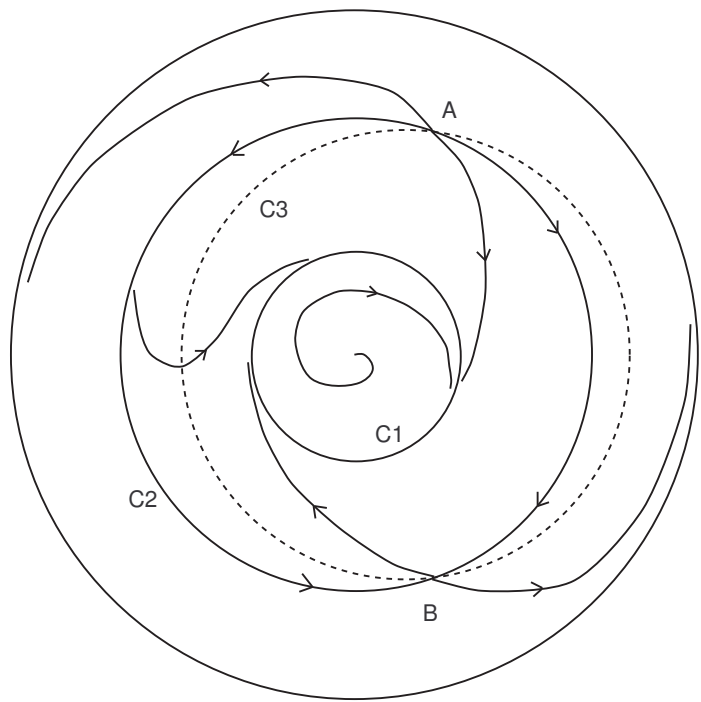,height=3.8in}}
\caption{Sketch of the global phase portrait corresponding to
Poincar\'{e}'s analysis of example P4. Circles $C_1$,$C_2$, and $C_3$,
which are referred to later in this paper, are labeled.  (The dashed line
for $C_3$ indicates that it is a nullcline and not a trajectory.) }
\end{figure}

To check the assertion, one applies Theorem 2 to the system.  By looking
at the the expression for P4 in polar coordinates, one finds that $I=5$
and $J=2$, hence $k=3$ and every point at infinity is a fixed point. There
is no limit cycle at infinity because every point on the equator is a
fixed point.  (If one had used the second convention, which was used by
Poincar\'e, one would have found trajectories crossing the equator.  With
this convention, the trajectories crossing the equator eliminate the
possibility of a cycle on the equator.  All reference to this second
convention will be dropped for all of the following examples.) The correct
global phase portrait for the system is presented in Figure~3.2.

\begin{figure}
\centerline{\psfig{file=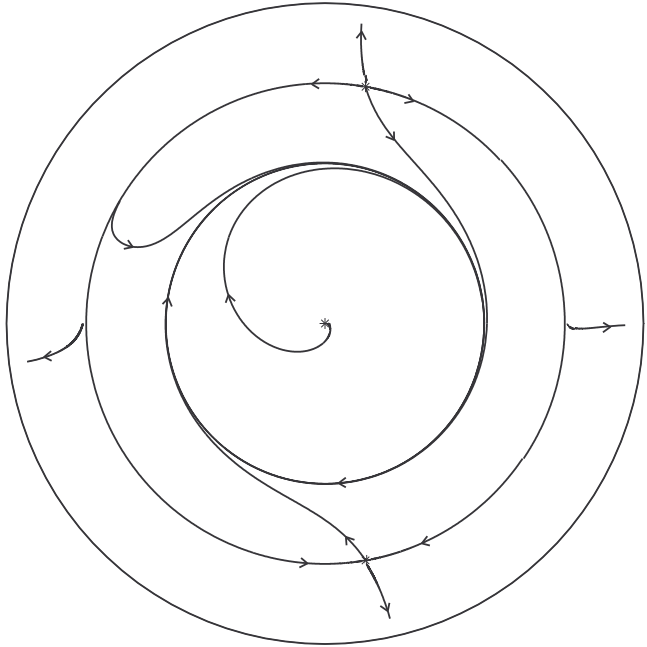,height=4.3in}}
\caption{The correct global phase portrait corresponding to Example P4.
	(Figures 3.2 and 4.1 were calculated numerically using a Runge Kutta
	method.)}
\end{figure}

\vspace{2 mm}
\newpage
\noindent{\bf Example P5}

Poincar\'{e}'s fifth example is not as clear of an extension of his third
and fourth examples as it does not quite fit the same ``genre'' of system.  
The main difference is that whereas P3 and P4 are more naturally
understood in a polar coordinate system, P5 is more easily understood in a
coordinate system based on lemniscates.  Example P5 is also more
complicated because it depends upon the parameter $c$.  He considers the
equation:
$$\frac{dx}{x(2x^2+2y^2+1)((x^2+y^2)^2+x^2-y^2-c)-y(2x^2+2y^2-1)}=$$
$$  \frac{dy}{y(2x^2+2y^2-1)((x^2+y^2)^2+x^2-y^2-c)+x(2x^2+2y^2+1)}.$$
Poincar\'{e} asserts that P5 has the following behavior, which depends on the 
parameter c.
\begin{enumerate}
\item For $c \le -1/4$, the system has:
\begin{itemize}
\item A saddle point at the origin.
\item Two unstable spirals at $(0,\pm1/\sqrt2)$.
\item An attracting limit cycle at infinity.
\end{itemize}
\item For $-1/4 < c < 0$, the system has:
\begin{itemize}
\item A saddle point at the origin.
\item Two stable spirals at $(0,\pm1/\sqrt2)$.
\item Two repelling limit cycles surrounding each of the two above spirals
	respectively.
\item An attracting limit cycle at infinity.
\end{itemize}
\item For $c=0$, the system has:
\begin{itemize}
\item A saddle point at the origin.
\item Two stable spirals at $(0,\pm1/\sqrt2)$.
\item Two homoclinic orbits which branch from the saddle and surround each of
	the spirals respectively.
\item An attracting limit cycle at infinity.
\end{itemize}
\item For $c>0$, the system has:
\begin{itemize}
\item A saddle point at the origin.
\item Two stable spirals at $(0,\pm1/\sqrt2)$.
\item One single repelling limit cycle centered at the origin which 
	surrounds all three fixed points.
\item An attracting limit cycle at infinity.
\end{itemize}
\end{enumerate}

The four qualitatively different phase portraits are sketched in figure 3.3.
\begin{figure}
\centerline{\psfig{file=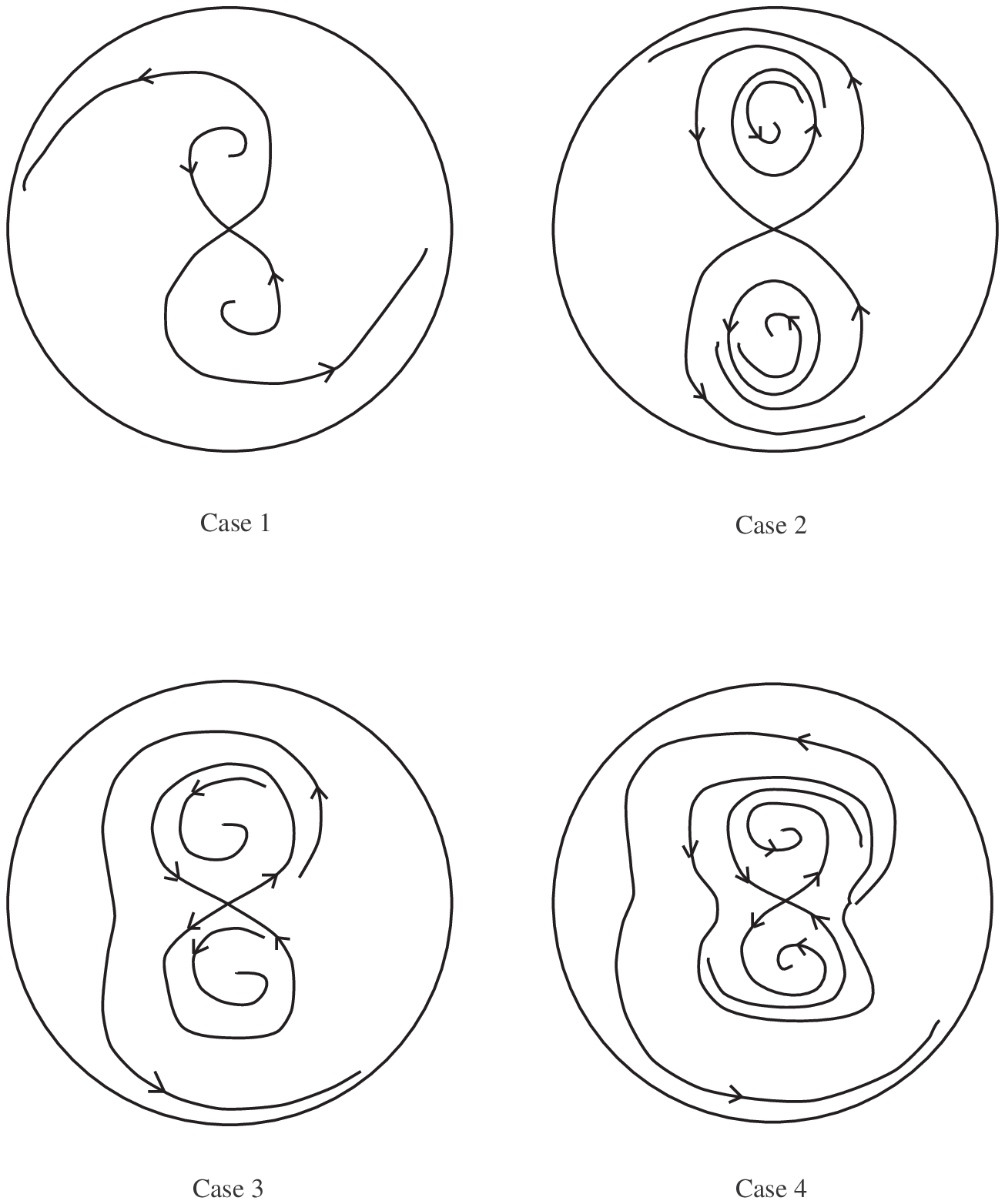,height=4.3in}}
\caption{Global phase portraits corresponding to Poincar\'e's analysis for 
the four different cases of P5.}
\end{figure}
As in the example P4, all of these features in the finite plane can be
easily verified. Because Poincar\'{e} does not include his analysis of the
behavior at infinity we check it here. Example P5 does not simplify in
polar coordinates, so it is easiest to use Theorem 1 to determine the
behavior at infinity. The highest order terms are of degree seven so we
have:
\begin{eqnarray*}
P_7 & = & 2X(X^2+Y^2)^3, \\
Q_7 & = & 2Y(X^2+Y^2)^3.
\end{eqnarray*}
So to find the fixed points at infinity, according to Theorem 1, one must 
solve the following system:
\begin{eqnarray*}
XQ_7-YP_7  =  2XY(X^2+Y^2)^3-2YX(X^2+Y^2)^3 & = & 0, \\
X^2+Y^2 & = & 1.
\end{eqnarray*}
The first equation is satisfied for all $X$,$Y$, so every point at
infinity is a fixed point.  Hence there can be no limit cycle at infinity.  
(This error was first discovered by Professor Lawrence Perko of Northern
Arizona University when he read a preliminary draft of this paper.)

\vspace{2 mm}

\noindent{\bf Examples from chapter nine}

Poincar\'e presents three additional examples in chapter nine
\cite[p.292-296]{PON} of his paper to illustrate the more general
techniques that he developed in chapter eight.  Poincar\'e claims that
each of these examples has a limit cycle at infinity.  However, one can
easily check his claim using Theorem 1 and Theorem 2, as was done above
for P3, P4, and P5, to see that each of these examples has fixed points at
every point on the equator. Because it is easy to check this, the
verification is left to the reader.  Furthermore, throughout this paper,
we focus on a discussion of examples P4 and P5.

The fact that these errors were so easily detected immediately makes one
wonder whether Poincar\'e used a similar form of analysis for finding
fixed points at infinity as was presented in Section 2.  When one looks 
at his first and second examples \cite[p. 274-278]{PON} one finds that he
uses the equations from Theorem 1 to find fixed points at infinity.
Because the calculations for finding fixed points at infinity for 
each of these examples are so simple (using Theorem 1 or Theorem 2) 
one is led to speculate as to the cause of Poincar\'e's error.

\setcounter{figure}{0}
\section{Error of Omission}
The facts that Poincar\'{e}'s errors in finding the fixed points at
infinity for his examples P4 and P5 were so easily detected, and that he
used the same type of criterion that is used here, suggest that, perhaps,
his analysis was of different equations than the ones printed in the
paper.  Were there to be an error of omission, terms which would lead to
an unsolvable system for finding the fixed points at infinity could be
missing which would not change the qualitative nature of his system on the
finite plane, where his analysis matched the printed equation.  If this
were the case, there most probably was an error of omission in the
production of his manuscript, not in the mathematics.  After lengthy
investigation, a system which could very plausibly become the system
printed in the paper as P4, were one to omit certain terms, becomes
convincing candidates for the intended fourth example.
	
In terms of the qualitative features predicted, Poincar\'{e}'s fourth
example appears to be an augmented version of his third example; however,
the equation listed is not an algebraically augmented version of the
equation in his third example. Algebraically, when in polar coordinates,
the $dr/dt$ equation is augmented with another nullcline at $r=3$.  To
create the new fixed points $A$ and $B$, the $d\theta /dt$ equation has a
new nullcline which crosses the $dr/dt$ nullcline at $r=3$.  However, the
$d\theta /dt$ equation is missing the factor $(r^2+1)$, which would not
create any nullclines, but which was necessary for P3 to have a limit
cycle at infinity (see calculation.)  A plausible attempt at finding
Poincar\'{e}'s actual equation is to include this term in the polar
version of example P4.  This attempt would fit with the intuition that P4
is merely an augmented version of P3.  (Example P5 is not considered ``in
the context of examples P3 and P4'' because the limit sets for P5 are
based on lemniscates instead of circles.)

\vspace{2 mm}

\noindent{\bf Example R4}

The system becomes:
\mathindent=20pt  
\begin{eqnarray*}
\frac{dr}{dt} & = & r(r^2-1)(r^2-9), \\
\frac{d\theta }{dt} & = & (r^2-2rcos(\theta ) - 8)(r^2+1).
\end{eqnarray*}
In Cartesian coordinates the system becomes:
\mathindent=20pt  
\begin{eqnarray*}
\frac{dx}{dt} & = & x(x^2+y^2-1)(x^2+y^2-9)-y(x^2+y^2-2x-8)(x^2+y^2+1), \\
\frac{dy}{dt} & = & y(x^2+y^2-1)(x^2+y^2-9)+x(x^2+y^2-2x-8)(x^2+y^2+1).
\end{eqnarray*}
This new version of Poincar\'{e}'s fourth example has all of the
qualitative features that his printed system has on the finite plane, and
has the additional feature of a limit cycle, instead of fixed points, at
infinity.  Hence this system exactly matches Poincar\'{e}'s analysis of
P4. The global phase portrait corresponding to example R4 is included in
figure 4.1. Verification that this system has the same qualitative
features as Poincar\'e claims for example P4 will not be presented here
because it is a result of Theorem 3 in the following section.  Because
this system algebraically follows from Poincar\'e's third example,
according to the hypothesis that Poincar\'e merely augmented his third
example to obtain his fourth example, this compels the author to believe
that this was Poincar\'e's intended system.

\begin{figure}
\centerline{\psfig{file=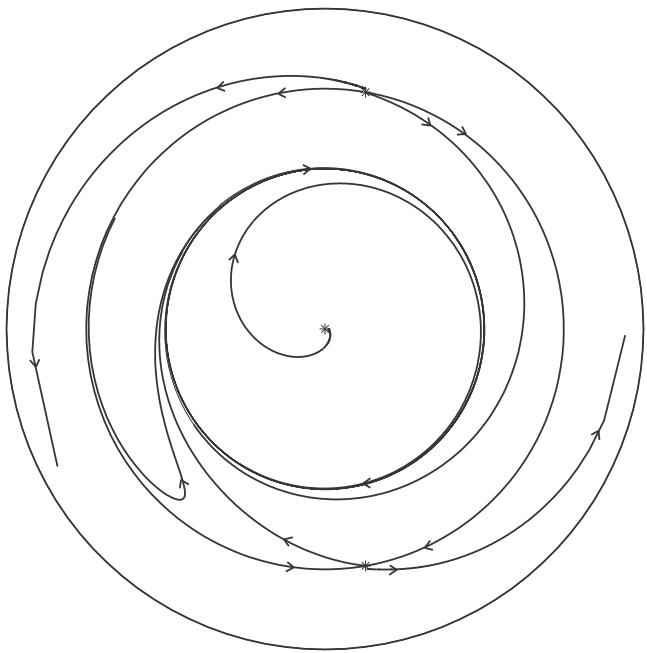,height=4.3in}}
\caption{The global phase portrait corresponding to example R4.}
\end{figure}

\section{A Method of Classification}

The results of the last section raise the question: how can one change the
algebraic statement of a differential equation and have the qualitative
aspects in the finite plane remain the same, while changing the behavior
at infinity? Rather than using {\it ad hoc} methods for changing an
equation, can there be a more systematic way of finding a certain
modification of an equation with the desired consequences?

To begin to answer these questions this section contains a systematic way
for finding which, among a certain class of changes, cause the equation
listed in Poincar\'{e}'s paper as P4 to be qualitatively equivalent to his
analysis listed for P4.  Based on the success of this technique with P4,
it is applied to P5 leading to a recommended modification of P5 which we
will call R5. These results not only serve as a demonstration of such
techniques, but also as a statement about the remarkable non-uniqueness of
a equations with this specific phase portrait.

This is the modification of P4:
\begin{theorem} {\bf (Classification Theorem for P4)}
Let $W(x,y)$ be a polynomial in $x$ and $y$ of degree $N$, and let
$W_N(x,y)$ be the homogeneous portion of $W$ which is of degree $N$.  If
\begin{enumerate}
\item $W(x,y) > 0$, for all $(x,y) \in \bf{R}^2$,
\item $N \geq 2$, and
\item $W_N(x,y)$ is positive definite.
\end{enumerate}
then:
the phase portrait of the system given by:
\begin{eqnarray*}
\frac{dx}{dt} & = & x(x^2+y^2-1)(x^2+y^2-9)-y(x^2+y^2-2x-8)W(x,y), \\
\frac{dy}{dt} & = & y(x^2+y^2-1)(x^2+y^2-9)+x(x^2+y^2-2x-8)W(x,y).
\end{eqnarray*}
gives a phase portrait qualitatively equivalent to the analysis Poincar\'{e}
gives in example P4.
\end{theorem}

First, one must verify that the system has the same number, location, and
type of fixed points.  To do this, first convert the system to polar
coordinates, obtaining:
\begin{eqnarray*}
\frac{dr}{dt} & = & r(r^2-1)(r^2-9), \\
\frac{d\theta}{dt} & = & (r^2-2rcos\theta -8)W(rcos \theta ,rsin \theta ).
\end{eqnarray*}
Because the new factor in the $d\theta /dt$ equation is never zero, it
does not create any new nullclines.  Hence, because there are no new
nullclines, only the fixed points that existed in the equation listed in
the paper ($A$,$B$, and $O$) exist in these modified examples.

To verify that these fixed points have the same local behavior as presented
by Poincar\'{e} one must calculate the trace, $\tau $, and determinant, $\Delta, $
of the Jacobian matrix of the system (in rectangular coordinates) 
symbolically at arbitrary $(x_0,y_0)$.
To simplify the calculation consider the following definitions:
\begin{eqnarray*}
C_1(x,y) & = & x^2+y^2-1, \\
C_2(x,y) & = & x^2+y^2-9, \\
C_3(x,y) & = & x^2+y^2-2x-8.
\end{eqnarray*}
With these definitions our system becomes:
\begin{eqnarray*}
\frac{dx}{dt} & = & xC_1(x,y)C_2(x,y)-yC_3(x,y)W(x,y), \\
\frac{dy}{dt} & = & yC_1(x,y)C_2(x,y)+xC_3(x,y)W(x,y).
\end{eqnarray*}
The expressions are calculated as follows (the arguments of circles $C_1$,
$C_2$, and $C_3$ are dropped):
\begin{eqnarray*} 
\tau & = & \frac{\partial P}{\partial x} + \frac{\partial Q}{\partial y}, \\
\Delta & = & 
 \frac{\partial P}{\partial x} \frac{\partial Q}{\partial y} - 
\frac{\partial P}{\partial y}  \frac{\partial Q}{\partial x}.
\end{eqnarray*} 
One obtains:
\begin{eqnarray*}
\tau & = & 2C_1C_2+2(x^2+y^2)(C_1+C_2)+2yW(x,y)+\\
 &  &  \hspace {.8in} C_3(x\frac{\partial W(x,y)}{\partial y}-y\frac{\partial W(x,y)}{\partial x}), \\
\vspace {.3in}
\Delta & = & \begin{array}{clcr}
 (C_1C_2+2x^2(C_1+C_2)-y(2x-2)W(x,y)-yC_3\frac{\partial W(x,y)}{\partial x})\cdot \\
 (C_1C_2+2y^2(C_1+C_2)+2xyW(x,y)+xC_3\frac{\partial W(x,y)}{\partial y}) - \\
 (2xy(C_1+C_2)-C_3W(x,y)-2y^2W(x,y)-yC_3\frac{\partial W(x,y)}{\partial y})\cdot \\
(2xy(C_1+C_2)+C_3W(x,y)+x(2x-2)W(x,y)+xC_3 \frac{ \partial W(x,y)}{\partial x}).
\end{array}
\end{eqnarray*}

When one evaluates these quantities at one of the fixed points, to check
the linearization, one is confronted with the unknown function $W(x,y)$.  
One must make restrictions on $W(x,y)$ so that the linearizations satisfy
Poincar\'{e}'s analysis for P4.

At fixed point $O$ all terms multiplied by $x$ or $y$ vanish, so one
obtains: $\tau = 18$, $\Delta=81+64W(0,0)^2$, and $\tau ^2-4\Delta =
-256W(0,0)^2$.  Because $\tau >0$ and because $\tau ^2-4\Delta <0$
(because $W(0,0)>0$) we have an unstable spiral.

At fixed point $B$ we have $\Delta = -144\sqrt{35} W(1/2,-\sqrt{35}/2)$.
Once again, because $W(x,y)>0$ for all $(x,y)$, we have $\Delta <0$,
hence a saddle.  

At fixed point $A$ we have $\tau = 144+\sqrt{35} W(1/2,\sqrt{35}/2)$ and
\newline $\Delta = 144\sqrt{35}W(1/2,\sqrt{35}/2).$ Because $W(x,y)>0$ we
have $\tau >0$.  Furthermore, $$\tau ^2-4\Delta =
(144-\sqrt{35}W(1/2,\sqrt{35}/2))^2 \geq 0.$$ Because $\tau^2-4\Delta \geq
0$ and $\tau > 0$, we find that fixed point $A$ is an unstable node.

Hence for any $W$ satisfying the conditions that $W(x,y)>0$ for all
$(x,y)$ then the fixed points of the modified system have the same
behavior as listed in Poincar\'{e}'s analysis of P4.

It is easy to verify that there is a stable limit cycle at $r=1.$  One 
must notice that $dr/dt=r(r^2-1)(r^2-9)=0$ for $r=1$, $dr/dt > 0$
for $0 < r < 1$, $dr/dt < 0$ for $1 < r < 3$, and that $d\theta/dt = 
(r^2-2rcos\theta -8)W(rcos\theta,rsin\theta) < 0$ for $0 < r < 2.$
Hence $r=1$ is a stable limit cycle.  

To establish the existence of the heteroclinic orbits, notice that for $r=3$,
 $dr/dt =0$. On this circle of radius three, to the left of $A$ and $B$,
$d\theta /dt > 0$, and to the right of $A$ and $B$, $d\theta /dt < 0$.  
The sign of $d\theta /dt$ remains the same in these regions as in the
printed equation because $W(x,y) > 0$ for all $(x,y)$ in $\bf{R}^2$.  
Hence there are trajectories which maintain $r=3$, connecting $A$ to $B$
both directions along the circle.

Finally we apply Theorem 2 to establish the existence of a limit cycle at 
infinity.  We can express $W(x,y)$ in polar coordinates in the following way:
$$W(rcos\theta,rsin\theta)=W_0(cos\theta,sin\theta)+rW_1(cos\theta,
	sin\theta)+\cdots+r^NW_N(cos\theta,sin\theta).$$
With $W_N(cos\theta,sin\theta) > 0$ for all $\theta$ because $W_N(x,y)$ is
positive definite.  Hence, using the notation of Theorem 2, $I=5$ and
$J=2+N$, where the degree of $W(x,y)$, $N \geq 2$.  Furthermore,
$\xi_J(\theta)=W_N(cos\theta,sin\theta) \neq 0$ for all $\theta$.  Hence,
there is a cycle at infinity. This cycle is a limit cycle because for $r
\geq 4$, $dr/dt > 0$.

This concludes the proof of Theorem 3.

Based on the success of the above method one may wish to see how well it
works for another system--we use P5.
\begin{theorem} {\bf (Classification Theorem for P5)}
Let $W(x,y)$ be a polynomial in $x$ and $y$ of degree $N$, and let
$W_N(x,y)$ be the homogeneous portion of $W$ which is of degree $N$.  If
\begin{enumerate}
\item $W(x,y) > 0$, for all $(x,y) \in \bf{R}^2$,
\item $N \geq 4$, and
\item $W_N(x,y)$ is positive definite,
\end{enumerate}
then: the phase portrait of the system given by:
\mathindent=0pt  
\begin{eqnarray*}
\frac{dx}{dt} & = & x(2x^2+2y^2+1)((x^2+y^2)^2+x^2-y^2-c)-y(2x^2+2y^2-1)W(x,y), \\
\frac{dy}{dt} & = & y(2x^2+2y^2-1)((x^2+y^2)^2+x^2-y^2-c)+x(2x^2+2y^2+1)W(x,y),
\end{eqnarray*}
gives a phase portrait qualitatively equivalent to the analysis Poincar\'{e}
gives in example P5.
\end{theorem}

We follow the same technique as in the previous theorem; first we verify 
that the modification results in the same fixed point behavior.  Clearly,
$W(x,y)$ is never zero, so it does not lead to any new fixed points.  Hence
we must only check that the existing fixed points have the same behavior.
To make the analysis easier, we make the following definitions:
\begin{eqnarray*}
A(x,y) & = & x(2x^2+2y^2+1), \\
B(x,y) & = & y(2x^2+2y^2-1), \\
C(x,y) & = & (x^2+y^2)^2+x^2-y^2-c.
\end{eqnarray*}
With these definitions the system becomes:
\begin{eqnarray*}
\frac{dx}{dt} & = & A(x,y)C(x,y)-B(x,y)W(x,y), \\
\frac{dy}{dt} & = & B(x,y)C(x,y)+A(x,y)W(x,y).
\end{eqnarray*}
At the fixed points $A(x,y)=B(x,y)=0$ and $x=0$ so in the expressions for
the trace and the determinant we have:
\begin{eqnarray*}
\tau & = & 8y^2C(x,y), \\
\Delta & = & (2y^2+1)(6y^2-1)((C(x,y))^2+(W(x,y))^2)
\end{eqnarray*}

At $(0,0)$ we have $\Delta=-(c^2+W(0,0)^2)$ which is negative for all
$(x,y)$ because $W(0,0)>0$.  Hence we have a saddle at the origin,
independent of the parameter value $c$ and the modification $W$.

At $(0,\pm 1/\sqrt2)$ we have $\tau=4C(0,\pm1/\sqrt2)=-4(1/4+c)$ and
\newline $\Delta= 4C(0,\pm1/\sqrt2)^2+4W(0,\pm1/\sqrt2)^2$.  Hence,
$\tau^2-4\Delta=-16W(0,\pm1/\sqrt2)^2 < 0$ resulting in spirals for any
choice of $W$.  The stability depends upon $\tau$; here $\tau > 0$ for $c
< -1/4$ and $\tau < 0$ for $c>-1/4$.  This matches Poincar\'{e}'s
assertion of unstable spirals for $c<-1/4$ and stable spirals for
$c>-1/4$.  Hence, as long as $W(x,y) > 0$ for all $(x,y)$ the behavior of
the fixed points is unchanged by $W$.

Next, we demonstrate that $C(x,y)=0$ is a conserved quantity.  This results
in the limit cycle(s) of cases 2 and 4 and the homoclinic orbits of case 3.
Consider the rate of change in $C(x,y)$ with respect to time:
\begin{eqnarray*}
\frac{\partial C(x,y)}{\partial t} & = & \frac{\partial C(x,y)}{\partial x}
  \frac{\partial x}{\partial t}+\frac{\partial C(x,y)}{\partial y}
  \frac{\partial y}{\partial t} \\
 & = & \frac{\partial C(x,y)}{\partial x}
  (AC-BW)+\frac{\partial C(x,y)}{\partial y}
  (BC+AW).
\end{eqnarray*}
Now one makes the observation that $\frac{\partial C}{\partial x}=2A$ and
$\frac{\partial C}{\partial y}=2B$.  So that: $$\frac{\partial C}{\partial
t}=2(A^2+B^2)C.$$ So, clearly $C(x,y)=0$ is a conserved quantity.  At this
point we see that the homoclinic orbits for case 4 have been verified
because for $c=0$ the algebraic curve $C(x,y)=0$ is a leminescate based at
the origin.  Further, for $C<0$, we have $\frac{\partial C}{\partial t}
<0$ and for $C>0$ we have $\frac{\partial C} {\partial t} >0$.  Based on
these results, one may find the necessary ``trapping regions'' to prove
the existence of the limit cycles using the Poincar\'e-Bendixson Theorem.

Finally, Poincar\'{e} asserted that for all of these values of $c$ there
is a a limit cycle at infinity.  To verify this, we determine whether
there are fixed points at infinity by using Theorem 1.  Using that $N \geq
4$, we obtain:

\begin{eqnarray*}
XQ_m-YP_m & = & 2(X^2+Y^2)^2W_N(X,Y)=0, \\
X^2+Y^2 & = & 1.
\end{eqnarray*}
which is clearly inconsistent because $W_N(X,Y)$ is positive definite, by
hypothesis.

Now that it has been verified that there are no fixed points at infinity,
to prove that there is a limit cycle at infinity it must be shown that
trajectories near infinity do, in fact, approach the cycle at infinity. In
this case we must apply slightly different techniques than Corollary 1
because $dr/dt$ is not single signed for all values of $\theta$.  
However, we only must show that the system's behavior near infinity does
approach the cycle at infinity.  To do this we note that $dC/dt =
2(A^2+B^2)C > 0$ for all points in the plane which are outside of the
limit cycle(s) or homoclinic orbits defined by the algebraic equation
$C(x,y) = 0$.  Based on this, one can rule out other limiting sets in the
proximity of infinity.  Because there are no other limiting sets within a
certain distance from the equator, the trajectories near the equator must
approach the cycle at the equator--the limit cycle at infinity.

This concludes proof of Theorem 4.

These theorems demonstrate the incredible algebraic non-uniqueness of the
equations which have the qualitative features outlined by Poincar\'{e} for
P4 and P5.  This non-uniqueness might be of concern in speculating which
modifications are the right ones. However, were one to speculate, without
insight into Poincar\'{e}'s other examples, as to which of these possible
modifications to P4 might be his equation, one might first try the most
simple of all $W$, $W(x,y)=(x^2+y^2+1)$.  This gives example R4 exactly.  
It is interesting that through choice by simplicity (application of
Occam's Razor) one obtains the same result which was found through
analysis of the context of the problem. This interesting philosophical
observation continues to persuade the author that R4 really was the system
intended by Poincar\'{e} for P4.

Based on all of the above success, the author would like to suggest an
adequate system to replace P5:

\vspace{2 mm}
\noindent{\bf Example R5}
\newline
The system becomes:
$$\frac{dx}{dt}  =  x(2x^2+2y^2+1)((x^2+y^2)^2+x^2-y^2-c)-y(2x^2+2y^2-1)(x^4+y^4+1), $$
$$\frac{dy}{dt}  =  y(2x^2+2y^2-1)((x^2+y^2)^2+x^2-y^2-c)+x(2x^2+2y^2+1)(x^4+y^4+1). $$
The fact that R5 matches Poincar\'{e}'s analysis of P5 is a clear application 
of Theorem 4.

\section{Discussion}

The reader may wonder why these problems were investigated. This section
briefly answers this question and attempts to bring closure to the
problem.

Verification of Poincar\'{e}'s assertions about his fourth example were
given as part of the take home project for a differential equations class.  
We were given the English translation as a guide \cite{DAV}.  One aspect
of the assignment was to attempt, through the use of a computer, to obtain
an actual global phase portrait. When this was done, the trajectories were
found to not cycle around near the equator, but to come in virtually
orthogonal to it! (See figure 3.2.)  It is by this means that the error
was first suspected.  Much speculation took place before all of the
analysis of P4 was done.

After a complete analysis of P4 came about, the error in example P5 was
discovered by L. Perko.  With example P5 in need of a modification and
with faith in the classification technique used for P4 (despite the fact
that it only considers one type of modification) the author set out to do
a similar analysis for P5.  It is by means of this type of analysis that
R5 was quickly obtained.

The errors found from chapter nine were found significantly later.  
Because the errors made with these equations are consistent with the
previous errors, little analysis has been done other than identifying the
errors.

The author finds it convincing that Poincar\'e could have made an error of
omission on his fourth example.  The relation between the third and fourth 
examples strongly supports this belief and the fact that of many plausible
modifications of Poincar\'e's fourth example, R4 is the simplest makes it 
convincing.  However, the fifth example does not fit within this genre of
system and it is not clear that the error made was one of omission.  The
further errors in chapter nine make it convincing that the remaining errors
where merely standard mathematical errors.

The author would like to conclude with a quotation of Poincar\'e: 

``How is an error possible in mathematics?  A sane mind should not be guilty 
of a logical fallacy, yet there are some very fine minds incapable of 
following mathematical demonstrations.  Need we add that mathematicians
themselves are not infallible?'' \cite{WEB}

\vspace{2 mm}

\noindent{\bf Acknowledgments}

I would like to thank Professor Jay P. Fillmore of UCSD for providing
insight, guidance, and most importantly inspiration to me for this
project.  If it weren't for his efforts, this paper would not exist.  I
would like to thank Professor Donald Smith of UCSD for providing many
helpful suggestions about this paper. I would also like to thank the
teaching assistant Greg Leibon for helping me accept that there actually
could be an error in Poincar\'{e}'s example.  I would like to thank
Professor Lawrence Perko of Northern Arizona University for his
suggestions and especially for finding the error in P5.  Finally, I would
like to thank both Professor John Guckenheimer and Professor Richard Rand
from Cornell University for the discussions that I had with each of them
regarding the analysis of differential equations at infinity.

\bibliography{p}
\bibliographystyle{plain}

\end{document}